\documentclass{article}
\usepackage{quotes}
\usepackage{algpseudocode}
\usepackage{algorithm}
\usepackage{rotating}
\usepackage{xcolor}
\usepackage{caption}

\usepackage{quotes}
\usepackage[utf8]{inputenc}
\usepackage{booktabs}
\usepackage{amsfonts}
\usepackage{epsfig}
\usepackage{epstopdf}
\usepackage{amsmath}
\usepackage{longtable}
\usepackage{graphicx}
\usepackage{amssymb}

\usepackage{algpseudocode}
\usepackage{lscape}
\usepackage{algorithm}
\usepackage{algorithmicx}

\newtheorem{theorem}{Theorem}

\newtheorem{proposition}[theorem]{Proposition}

\newenvironment{proof}[1][Proof]{\noindent\textbf{#1.} }{\ \rule{0.5em}{0.5em}}

\newcommand{\be}{\begin{enumerate}}
	\newcommand{\ee}{\end{enumerate}}
\newcommand{\bc}{\begin{center}}
	\newcommand{\ec}{\end{center}}
\newcommand{\bi}{\begin{itemize}}
	\newcommand{\ei}{\end{itemize}}

\newcommand{\TSP}{\textsc{TSP}}
\newcommand{\PTP}{\textsc{PTP}}
\newcommand{\PPTP}{\textsc{PPTP}}

\usepackage{subcaption}
\usepackage{multirow}
\title{The Probabilistic Profitable Tour Problem\\ under a specific graph structure}
\author{Enrico Angelelli$^{(1)}$, Renata Mansini$^{(2)}$, Romeo Rizzi$^{(3)}$\\ 
{\small	$(1)$ Department of Economics and Management, University of Brescia, Italy; 
	{\tt enrico.angelelli@unibs.it}}\\
{\small	$(2)$ Department of Information Engineering, University of Brescia, Italy;
	{\tt renata.mansini@unibs.it}}\\
{\small	$(3)$ {Department of Computer Science, University of Verona, Italy};
	{\tt romeo.rizzi@univr.it}}\\
}

\begin{document}

\maketitle		
		
		\begin{abstract}
			Among the most important variants of the traveling salesman problem (\TSP)
			are those relaxing the constraint that every {\it locus} 
			should necessarily get visited,
			rather taking into account a revenue (prize) for visiting customers.
			In the Profitable Tour Problem (\PTP), we seek for a tour visiting a subset of customers
			while maximizing \emph{net gain} (profit) as difference between total revenue collected from visited customers and incurred traveling costs.
			The metric \TSP\ can be modeled as a \PTP\ with large revenues.
			As such, PTP is well-known to be NP-hard and also APX-hardness follows.
			Nevertheless, \PTP\ is solvable in polynomial time on particular graph structures like lines, trees and circles.
			Following recent emphasis on robust optimization,
			and motivated by current flourishing of retail delivery services,
			we study the Probabilistic Profitable Tour Problem (\PPTP), the generalization of \PTP\ where customers will show up with a known probability,
			in their respective loci,
			only after the tour has been planned.
			Here, the selection of customers has to be made \emph{a priori}, before knowing if a customer will actually submit his request or will not.
			While the tour has to be designed 
			without this knowledge, revenues will only be collected  from customers who will require the service. The objective is to maximize the expected net gain obtained by visiting only the customers that show up.
			We provide a polynomial time algorithm
			computing and characterizing the space of optimal
			solutions for the special case of the \PPTP\ where customers are distributed on a line.
			
		\end{abstract}
		

\section{Introduction}

Several variants of the traveling salesman problem, taking into account revenues (profits) for visiting
customers, have been studied in the literature: with single  or multiple vehicles, with and without time windows,
with precedence (Hanafi et al.\ \cite{HMZ20})
or other side constraints (see the surveys by Feillet et al.\ \cite{FDG05} and, more recently, by Gunawan et al\ \cite{GLV16}).

In the PTP,  each customer with a known location  is associated with 
a positive revenue (prize) and the problem looks for
a tour
originating at a depot and
visiting
a subset of the customers so to maximize the difference between collected revenues and total traveling costs (Feillet et al.\ \cite{FDG05}).

In the literature,
PTP is sometimes formulated as a minimization problem trying to find a tour that minimizes the sum of  cost and missed revenue for not visiting customers.  The two formulations complement each others and  are equivalent as optimization problems (Johnson et al.\ \cite{JMP00}).
A great interest has been devoted also to the capacitated variant of the PTP,  since it can be seen as a special case of the elementary shortest path problem  with resource constraints (see Jepsen et al.\ \cite{JPSP14} where the authors introduce  a branch-and-cut and new valid inequalities for the problem). Recently, Angelelli et al.\ \cite{ABST14} study the deterministic
PTP and  other variants of the Traveling Salesman problem with profits considering  special structures of the underlying graph (line, cycle, star rooted at the depot, tree rooted at the depot). The authors provide computational complexity  and approximation results for all studied problems generalizing them to the case with positive service times  associated with  customers.
As far as the profitable tour problem on a line (a path) is concerned, the authors show that it is solvable in linear time under both service settings.

Frequently, in many real application contexts, the customers will
require a service with a known probability.
In such a case, the decision maker has to decide a strategy  to construct an a priori solution, i.e.\ he/she has to select  which  customers to visit before knowing who among them will submit a  request. Only those customers,  among the  ones selected in advance,  who will submit a request,  will also be served.
The objective of the problem is to select a subset of customers so to maximize the expected profit measured as difference between  revenue and traveled distance.  We call this variant of the PTP, the Probabilistic Profitable Tour Problem (PPTP).

In contrast to the deterministic variants, a very few contributions can be found on probabilistic TSP with profits. Angelelli et al.\  \cite{AAFV17}
provide a linear integer stochastic formulation of the Orienteering Problem and develop both a branch-and-cut approach and different matheuristic methods. To the best of our knowledge, no
contributions can be found on the PTP under uncertainty.

In this paper, we will analyze the probabilistic PTP where customers are distributed on  a line.  Starting from main results provided in Angelelli et al.\ \cite{ABST14}, our  theoretical question is 
if the probabilistic PTP on a line remains solvable in polynomial time
or becomes NP-hard.
Our main contribution is the complexity of  probabilistic PTP on a line and the characterization of its optimal solutions space.
In particular, we show that the  problem  can be solved in  $O(n^3)$ time, where $n$ is the number of customers. 
Although the main contribution of our work remains theoretical, one can figure out real contexts where  the problem  is likely to find application: the road network topology is typical of  a mountain valley where customers are all located on a main road, whereas the a priori optimization is typical of a 2-stage decision process where some decisions are made at the first stage (selection of the  customers for whom service is guaranteed) and a recourse action takes place  at the second stage (maximize the expected value of the total collected revenues minus the traveling costs related to the subset of a priori selected customers that actually showed up by requiring the service).

The paper is organized as follows.
In Section \ref{sec:properties}, the main properties of the problem and its optimal solution are discussed for a given value of the 
unitary cost per  traveled distance, whereas in Section \ref{sec:dependency-c}  the dependency on such a parameter 
is analyzed and the space of optimal
solutions is characterized. Finally, the solution algorithm is presented in Section \ref{sec:algorithm} and its computational complexity discussed.
Concluding remarks are provided in Section \ref{sec:conclusion}.

\section{Problem definition and main properties}\label{sec:properties}

Formally, the PPTP can be defined on a directed graph
$G=(V,A)$  with $V= \{v_0\} \cup N$
where $v_0$ is the depot,
$N=\{v_1,...,v_n\}$ is the set of potential customers and $A$ denotes
the arc set.
To each customer $v_i$, a revenue (prize) $p_{i}>0$   and a probability $\pi_{i}\in(0,1]$
are assigned.
The depot $v_0$ is located at the origin of the semi-line in the point $x_0=0$. Each customer $v_i$ is positioned in  $x_i>0$.
Each arc $(i,j)$ is weighted by
a cost $|x_j-x_i|c$ where
$c \geq 0$ is a unitary cost per  traveled distance.
Without loss of generality, we assume that:
\begin{itemize}
	\item $0<x_{i}\leq x_{i+1}$
	for $i=1,...n-1$;
	\item $x_{i}=x_{i+1}\Rightarrow\pi_{i}%
	\leq\pi_{i+1}$ (if two customers are positioned at the same location, they are sorted in non decreasing  order of their probability)
	
	\item if two customers share the same location and have equal probability, the  tie is broken  randomly.
\end{itemize}

Let  $S\subseteq N$ the set of customers selected a priori by the decision maker and let $X$ be the subset of customers that will submit a request. The set of customer that will be actually served by the vehicle is $S \cap X.$ We indicate as $R(S)$
the corresponding expected revenue:

\begin{equation}
R(S)=\sum_{v_{i}\in S}\pi_{i}p_{i}.
\end{equation}

	Observe that the minimum length route  to serve all customers in $S \cap X$ corresponds  to reach  the farthest customer while serving all the others on the way back to the depot. This means that in the case $v_i$ is the farthest customer in $S \cap X$, the total traveled distance is $2x_i$, that, for simplicity of notation, we denote as $l_i$.
Thus, according to the assumption $0<x_{i}\leq x_{i+1}$,  the
expected traveled distance can be computed as:

\begin{equation}\label{eq:distance}
L(S)=\sum_{v_{i}\in S}\left[  l_{i}\pi_{i}%
{\displaystyle\prod_{{v_j}\in S,\ j>i}}
(1-\pi_{j})\right]  .
\end{equation}

Given the cost $c \geq 0$ per unit of traveled length, the expected cost $C(S,c)$ and the expected profit $G(S,c)$ will be computed as follows:
\begin{equation}
C(S,c)=cL(S),
\end{equation}
\begin{equation}
G(S,c)=R(S)-cL(S).
\end{equation}

The PPTP on a line can thus be formulated as follows:
\begin{equation}
PPTP(N,c)=\max_{S \subseteq N}G(S,c) \label{POP-def}
\end{equation}
We call a solution of (\ref{POP-def}) an \emph{optimal set}. The optimal set may depend on $c$, but by now we assume $c$ as given.

Solving problem (\ref{POP-def}) may appear a difficult task as the number of potential options is exponential with respect to $n$ for every chosen $c$.
However,  in Section \ref{sec:algorithm}, we will show an approach that, looking at $PPTP(N,c)$ as a function of $c$, manages to build in polynomial time a family of optimal sets covering all possible values of $c\in[0,+\infty)$.

In this section, we start showing some basic properties of the problem for a fixed value of $c$ and conclude that, even though we may have several optimal sets for a given $c$, there is only one \emph{maximal optimal set} and only one \emph{minimal optimal set}, whereas any other optimal set is a subset of the former and a superset of the latter.

\medskip

Given a non void subset of customers $ S\subseteq N,$ we refer to $v^{S} \in S$ as the customer with the highest index in $S$ (the farthest customer from the depot and with highest probability). Moreover, we  indicate as  $S^{\prime}$ the set
$S\backslash\{v^{S}\},$  whereas  $l^{S},$ $p^{S}$ and $\pi^{S}$ are distance, score and probability of $v^S$.

\begin{proposition}\label{PROP:RICORSIONE_1}
	Given a non void set of customers $S$, the following recursive formulas hold: (See proof in appendix)
	
	\begin{equation}
	R(S)=R(S^{\prime})+\pi^{S}p^{S} \label{RICORSIONE_R}%
	\end{equation}
	\begin{equation}
	L(S)=(1-\pi^{S})L(S^{\prime})+\pi^{S}l^{S} \label{RICORSIONE_L}%
	\end{equation}
	\begin{equation}
	C(S,c)=C(S^{\prime},c(1-\pi^{S}))+\pi^{S}l^{S}c \label{RICORSIONE_C}%
	\end{equation}
	\begin{equation}
	G(S,c)=G(S^{\prime},(1-\pi^{S})c)+\pi^{S}(p^{S}-l^{S}c) \label{RICORSIONE_G}%
	\end{equation}
	
\end{proposition}

\begin{proposition}\label{PROP_SUBSET_LENGTH_INEQUALITY}
	Given any two subsets of customers  $S_1$ and $S_2$ such that $S_1 \subseteq S_2$, then $L(S_1)\le L(S_2)$. (See proof in appendix)
\end{proposition}

\begin{proposition}\label{PROP_INSIEMI}
	Given any two subsets of customers  $S_1$ and $S_2$, the following results hold: (See proof in appendix)
	\begin{equation}\label{PROP_INSIEMI_R}
	R(S_1\cup S_2) = R(S_1)+R(S_2) - R(S_1\cap S_2),
	\end{equation}
	\begin{equation}\label{PROP_INSIEMI_L}
	L(S_1\cup S_2) \leq L(S_1)+L(S_2) - L(S_1\cap S_2).
	\end{equation}
\end{proposition}

\begin{proposition}\label{PROP_INSIEMI_3}
	Given any two subsets of customers  $S_1$ and $S_2$, we have
	\begin{equation}
	G(S_1,c)+G(S_2,c) \le G(S_1\cup S_2,c)+G(S_1\cap S_2,c).
	\end{equation}
\end{proposition}

\begin{proof}
	From Proposition \ref{PROP_INSIEMI} we get
	\begin{align*}
	G(S_1\cup S_2,c)  &  =R(S_1\cup S_2)-cL(S_1\cup S_2)\\
	&  \geq\left[  R(S_1)+R(S_2)-R(S_1\cap S_2)\right]  -c[L(S_1)+L(S_2)-L(S_1\cap S_2)]\\
	&  \geq\lbrack R(S_1)-cL(S_1)]+[R(S_2)-cL(S_2)]-[R(S_1\cap S_2)-cL(S_1\cap S_2)]\\
	&  \geq G(S_1,c)+G(S_2,c)-G(S_1\cap S_2,c).
	\end{align*}
\end{proof}

\begin{proposition}\label{PROP_INEQUNINT}
	Given any two subsets of customers  $S_1$ and $S_2$ such that $G(S_1,c)=G(S_2,c)$ for a given unitary cost
	$c,$ then
	\begin{equation}
	G(S_1,c)\leq\max(G(S_1\cup S_2,c),G(S_1\cap S_2,c)).
	\end{equation}
\end{proposition}

\begin{proof}
	By contradiction,  if $\max(G(S_1\cup S_2,c),G(S_1\cap S_2,c))<G(S_1,c)$ then we get
	$G(S_1\cup S_2,c)+G(S_1\cap S_2,c)<G(S_1,c)+G(S_2,c)$ in contrast with
	Proposition \ref{PROP_INSIEMI_3}.
\end{proof}

\begin{proposition}
	\label{PROP_UNINT} If two subset of customers $S_1$ and $S_2$ are optimal sets for a given unitary cost $c$, then $S_1\cup S_2$ and $S_1\cap S_2$ are also optimal sets.
\end{proposition}

\begin{proof}
	From Proposition \ref{PROP_INEQUNINT} we know that $G(S_1,c)\leq\max(G(S_1\cup
	S_2,c),G(S_1\cap S_2,c)),$ but for the optimality of  $S_1$ and $S_2$ we have
	$G(S_1,c)=\max(G(S_1\cup S_2,c),G(S_1\cap S_2,c))$, which proves that at least one between  $S_1\cap S_2$ and $S_1\cup S_2$ is optimal. We now show that both are optimal.
	
	From Proposition \ref{PROP_INSIEMI_3} we know that
	\[
	G(S_1,c)+G(S_2,c)\leq\max(G(S_1\cup S_2,c),G(S_1\cap S_2,c))+\min(G(S_1\cup S_2,c),G(S_1\cap
	S_2,c)),
	\]
	but from equality  $G(S_1,c)=\max(G(S_1\cup S_2,c),G(S_1\cap S_2,c))$ it follows that
	\[
	G(S_2,c)\leq\min(G(S_1\cup S_2,c),G(S_1\cap S_2,c))
	\]
	and thus, from the optimality of $S_2$, the optimality of both $S_1\cap S_2$ and $S_1\cup S_2$ follows.
\end{proof}

\begin{proposition}\label{PROP_UNIQUE_MAX_MIN}
	For each unitary cost  $c$, there exists a unique minimal optimal set and a unique maximal optimal set. Moreover, any optimal set is a subset of the maximal one and a superset of the minimal one.
\end{proposition}

\begin{proof}
	By contradiction, let $S_1$ and $S_2$ be two distinct maximal optimal sets, we know from Proposition \ref{PROP_UNINT} that their proper superset $S_1\cup S_2$ is optimal, thus neither $S_1$ or $S_2$ can be maximal optimal sets.
	Moreover, if we had an optimal set $S$ which is not a subset of the maximal one $\bar S$ , then $S\cup\bar S$ would be optimal and a proper superset of $\bar S$, in contrast with maximality of $\bar S$.
	
	We use the same argument to show that the minimal optimal set is unique and any optimal set is a superset of the minimal one.
\end{proof}

\section{Dependency from the unitary cost}\label{sec:dependency-c}
In this section, we study the dependency of $PPTP(N,c)$ from the unitary cost $c$ and characterize the property of optimal sets accordingly. We end this section with a hint on how to build a description of $PPTP(N,c)$ in term of values and optimal sets.

In the following, we call {\it characteristic function} of a set of customers $S$ its expected revenue
$G(S,c)=R(S)-cL(S)$ seen as a function of the unitary cost $c.$
Note that the characteristic function of a set  $S$ can be graphically represented in the Cartesian plane $G(S,c)$ as a line with non positive slope.
In particular, the slope is  negative for each set  $S\neq\emptyset$ and null for
$S=\emptyset;$ in such a case  $G(\emptyset,c)=0$ for all  $c.$

\begin{proposition}
	\label{PROP_POP_C}  $PPTP(N,c)$ as a function of $c$ holds the following characteristics:
	
	\begin{enumerate}
		\item \label{PROP_POP_C_1}$PPTP(N,0)=R(N)=\sum_{v_{i}\in N}\pi_{i}p_i$
		
		\item \label{PROP_POP_C_2} There exists a finite positive value $\widetilde{c}$ such that
		\[
		\left\{
		\begin{array}
		[l]{c}%
		PPTP(N,c)>0\text{   if  }c<\widetilde{c}\\
		PPTP(N,c)=0\text{ otherwise} 
		\end{array}
		\right.
		\]
		
		\item \label{PROP_POP_C_3}In the interval $[0,\widetilde{c}]$ function $PPTP(N,c)$ is strictly decreasing, piece-wise linear and convex.  In $[0,\infty)$ function $PPTP(N,c)$ is non-increasing  piece-wise linear and convex.
	\end{enumerate}
\end{proposition}

\begin{proof}
	Let's analyze each point separately:
	
	\begin{enumerate}
		\item For  $c=0$ we get  $PPTP(N,0)=\max_{S\in P(N)}\left(  R(S)-0\cdot
		L(S)\right)  =\max_{S\in P(N)}\left(  R(S)\right)  =R(N).$
		
		\item For each non void set  $S\subseteq N$ we get $G(S,c)\leq0$ for
		$c\geq c_{S}=R(S)/L(S).$ Let us indicate  $\widetilde{c}=\max_{S\subseteq N}\{c_{S}\}$. Then we get that for each $c<\widetilde{c}$, inequality $PPTP(N,c)>0$ holds since there exists at least one set $S\subseteq N$ for which
		$G(S,c)>0;$ when $c\geq\widetilde{c}$ then we get $G(S,c)\leq0$ for each set $S\subseteq N$ and thus the void set corresponds to the optimal solution since
		$G(\emptyset,c)=0$.

		\item For $c\in\lbrack0,\widetilde{c}]$ function $PPTP(N,c)=\max_{S\in
			P(N)}G(S,c)$ is the envelop of a finite set of linear functions strictly decreasing. For $c\geq\widetilde{c}$, function
		$PPTP(N,c)$ becomes constant. In the interval $[0,\infty)$, it only loses monotonicity. 
	\end{enumerate}
\end{proof}

The function $PPTP(N,c)$ is thus described by a finite number of linear pieces.
Each linear piece is the characteristic function of an optimal set
within the corresponding range. To determine the actual shape of $PPTP(N,c)$ can, in principle, be hard as it is the outcome of a number,  exponential in $n$, of lines. However, we will show that this task can be accomplished in polynomial time $O(n^3)$. Let us start with some properties of function $PPTP(N,c)$.

\begin{proposition}\label{PROP_CORNER_POINT}
	If for some $\bar{c}$ there are two distinct optimal sets, then the maximal optimal set is optimal on $(\bar{c}-\varepsilon,\bar{c}]$ and the minimal optimal set is optimal on $[\bar{c},\bar{c}+\varepsilon)$ for some $\varepsilon>0$. In particular, $\bar{c}$ is a corner point of PPTP(N,c).
\end{proposition}

\begin{proof}
	Let us consider the minimal and maximal optimal sets in $\bar{c}$, and call them $A$ and $B$, respectively. By Proposition \ref{PROP_UNIQUE_MAX_MIN} we know that they are unique distinct and that $A\subset B$.
	It is easy to see that $R(A)<R(B)$ and $L(A)<L(B)$ (Proposition \ref{PROP_SUBSET_LENGTH_INEQUALITY} so that, being $G(A,\bar{c})=G(B,\bar{c})$ it must necessarily be $G(A,c)<G(B,c)$ for $c<\bar{c}$ and viceversa $G(A,c)>G(B,c)$ for $c>\bar{c}$. Easy to see that any other intermediate optimal set in $\bar{c}$ is dominated by the maximal optimal set for $c<\bar{c}$ and by the minimal optimal set for $c>\bar{c}$.
	
	Furthermore, if for a fixed $\varepsilon>0$ there were some distinct optimal sets in interval $(\bar{c}-\varepsilon,\bar{c}+\varepsilon)$ other than $A$ and $B$, then we can repeatedly halve the value of $\varepsilon$; the process must come to an end as we have only a finite number of potential optimal sets.
	
	Point $\bar{c}$ is a corner point because function $PPTP(N,c)$ takes different slopes around $\bar{c}$.
\end{proof}

\begin{proposition}\label{PROP_RANGE_OPTSET}
	If $(c_1,c_2)$ is an interval such that $PPTP(N,c)$ is linear (no corner points in the interval), then there is only one optimal set for all $c\in(c_1,c_2)$.
\end{proposition}

\begin{proof}
	If two distinct set of customers are optimal for some $c\in(c_1,c_2)$ then $c$ is a corner point.
\end{proof}

\begin{proposition}\label{PROP_N_CORNERS}
	There are at most $n$ corner points in function PPTP(N,c).
\end{proposition}

\begin{proof}
	According to Proposition \ref{PROP_CORNER_POINT}, at each corner point the optimal set loose some customers, and since we have $n$ customers, we can have $n$ corner points at most.
\end{proof}

\bigskip

Resuming, we showed that:
\begin{enumerate}
	\item In the corner points of function $PPTP(N,c)$, we have two optimal solutions defined by the linear pieces  belonging to the envelop, of which one is subset of the other. More precisely, the piece belonging to higher values of  $c$ is characterized by the minimal optimal  set which is a subset of the maximal optimal set characterized by the piece associated with lower values of  $c$ (Proposition \ref{PROP_CORNER_POINT});
	\item In each corner point $\bar{c}$ we may have other optimal sets, but all of them are dominated by the maximal optimal set for $c<\bar{c}$ and by the minimal optimal set for $c>\bar{c}$; each of these intermediate sets is a subset of the maximal one and a superset of the minimal one  (Proposition \ref{PROP_CORNER_POINT});
	\item Each linear piece of function  $PPTP(N,c)$ is characterized by one and only one optimal maximal solution (Proposition \ref{PROP_RANGE_OPTSET});
	\item given two consecutive corner points $c_1$ and $c_2$, and the line segment representing PPTP(N,c) for $c\in[c:1,c_2]$, the corresponding optimal set is the minimal optimal set for $c=c_1$, the only optimal set for $c\in (c_1,c_2)$ and the maximal optimal set in $c=c_2$ (Propositions \ref{PROP_CORNER_POINT} and \ref{PROP_RANGE_OPTSET});
	\item The function $PPTP(N,c)$ contains at most $n$ corner points and $n+1$ linear pieces (Proposition \ref{PROP_N_CORNERS});
	\item The function $PPTP(N,c)$ is defined for all $c>0$ by at most $n+1$ optimal sets.
\end{enumerate}
Next, we show that optimal sets of function $PPTP(N,c),$ (and the corresponding corner points) can be computed in polynomial time with respect to the size $n$ of the instance.

\bigskip

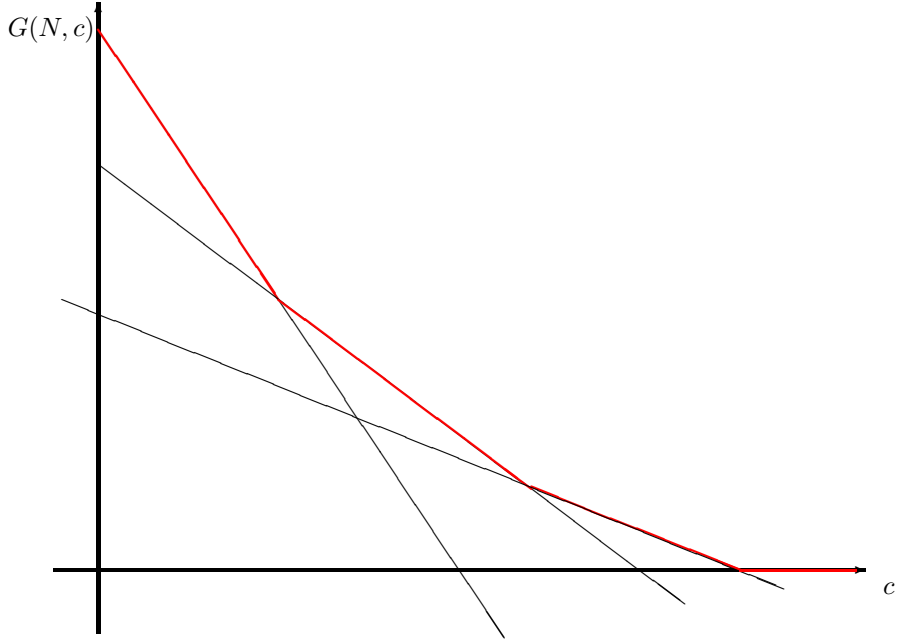
\begin{figure}[h]
	\begin{minipage}[t]{12cm}
		\setlength{\unitlength}{0.1\textwidth} \linethickness{1.1pt}
		\begin{picture}(7.00, 7.00)
		\put(0.5,1){\vector(1,0){9}} \put(1,0.3){\vector(0,1){7}}
		
		\put(-0.1,7){\makebox(0,0)[l]{{ $G(N,c)$}}}
		\put(1,7){\line(2,-3){4.5}}
		{\color{red}	
			\thicklines	\thicklines \thicklines\put(1,7){\line(2,-3){2}}}
		
		\put(1,5.5){\line(4,-3){6.5}}
		{\color{red}	\linethickness{1.5pt} \thicklines	\thicklines \thicklines\put(3,4){\line(4,-3){2.8}}}
		
		
		{\color{red}  \thicklines
			\thicklines \thicklines\put(5.7,1.928){\line(5,-2){2.35}}}
		
		\put(0.5,4){\line(5,-2){8}}
		{\color{red}  \thicklines
			\thicklines \thicklines \put(8,1){\line(1,0){1.3}}}

		\put(9.5,0.8){\makebox(0,0)[l]{{ $c$}}}
		\end{picture}
	\end{minipage}
	\caption{\label{fig:PPTP} The piece-wise linear and convex shape (in red) of function $PPTP(N,c)$: a simple case with $n=3.$  }
	
\end{figure}

By defining  $N^{(k)}$ as the set of the first  $k$ customers closest to the depot for $k=0,1,\ldots,n$,
we can also write:
\begin{equation}
PPTP(N^{(k)},c)=\max_{S\in \mathcal{P}(N^{(k)})}G(S,c),\ \ \ \ k=0,...,n,
\end{equation}
where in the particular case $k=0$ we get an empty optimal set with $PPTP(N^{(0)},c)=0$ for all $c\geq 0$. 
Our idea is to iteratively build function $PPTP(N^{(k+1)},c)$ from $PPTP(N^{(k)},c)$, starting with $k=0$ up to $k=n-1$.

Let us indicate with $\mathcal S^{(k)}$ the family of optimal sets defining function $PPTP(N^{(k)},c)$; We also indicate with $S^{(k)}_c\in \mathcal S^{(k)}$ the optimal set for a given unit cost $c$, that is
\begin{equation}
PPTP(N^{(k)},c)=\max_{S\in \mathcal S^{(k)}}(G(S,c))=G(S^{(k)}_c,c).
\end{equation}
Finally, we indicate with $\mathcal E^{(k)}=\{S\cup\{v_{k+1}\}\ |\ S\in \mathcal S^{(k)}\}$ the family of extended sets obtained from sets in $\mathcal S^{(k)}$ by adding a new customer $v_{k+1}$. We recall that according to Proposition \ref{PROP_N_CORNERS}, family $\mathcal S^{(k)}$ contains at most $k+1$ elements.

Next, we show that the family of optimal sets $\mathcal S^{(k+1)}$ for $PPTP(N^{(k+1)},c)$ can be extracted from $\mathcal S^{(k)}\cup \mathcal E^{(k)}$ which is the fundamental property the iterative step of the building process lays on.

\begin{proposition}
	\label{PROP_ESTENSIONE_K_K1} The family of optimal sets $\mathcal S^{(k+1)}$ for function $PPTP(N^{(k+1)},c)$ is a subset of $\mathcal S^{(k)}\cup \mathcal E^{(k)}$ with cardinality at most $(k+1)+1$.
\end{proposition}

\begin{proof}
	Let $S_1\subseteq N^{(k+1)}$ be a set of customers. We show that $G(S_1,c)\leq \max_{S \in \mathcal S^{(k)}\cup \mathcal E^{(k)}}(G(S,c))$ for all $c>0$. We proceed by cases.
	
	a) If $v_{k+1}\notin S_1,$ then by construction we have
	$$
	G(S_1,c)\leq PTPP(N^{(k)},c)=\max_{S\in \mathcal S^{(k)}}(G(S,c))\leq \max_{S\in \mathcal S^{(k)}\cup\mathcal E^{(k)}}(G(S,c)).
	$$
	
	b) If $v_{k+1}\in S_1,$ then let us define  $S_1^{\prime}=S_1\backslash\{v_{k+1}\}\subseteq
	N^{(k)}.$ By using (\ref{RICORSIONE_G}), we get:
	\begin{align*}
	G(S_1,c) &  =G(S_1^{\prime},(1-\pi_{k+1})c)+\pi_{k+1}(p_{k+1}-cl_{k+1})\\
	&  \leq\max_{S\in \mathcal S^{(k)}}(G(S,(1-\pi_{k+1})c))+\pi_{k+1}(p_{k+1}-cl_{k+1}),
	\end{align*}
	now let us indicate with $\widetilde{S}_c\in \mathcal S^{(k)}$ the optimal set of problem $\max_{S\in \mathcal S^{(k)}}(G(S,(1-\pi_{k+1})c))$ for a given $c$ and observe that $\widetilde{S}_c\cup \{v_{k+1}\}\in \mathcal E^{(k)}$; Thus, we have
	\begin{align*}
	G(S_1,c) &  \leq G(\widetilde{S}_c,(1-\pi_{k+1})c)+\pi_{k+1}(p_{k+1}-cl_{k+1})\\
	&  = G(\widetilde{S}_c\cup \{v_{k+1}\},c))\leq \max_{S\in \mathcal E^{(k)}}(G(S,c))\leq \max_{S\in S^{(k)}\cup E^{(k)}}(G(S,c))
	\end{align*}
	
	Thus, $\mathcal S^{(k)}\cup \mathcal E^{(k)}$ 
	is enough to determine function  $PPTP(N^{(k+1)},c).$ It is worth noticing that, according to Proposition \ref{PROP_N_CORNERS}, the family of customer sets
	$S^{(k)}\cup E^{(k)}$ can be reduced to have a cardinality no larger  than
	$k+2$.
\end{proof}

\section{Solution algorithm and computational complexity}\label{sec:algorithm}
This section is devoted to the presentation of the procedure to compute function
$PPTP(N,c).$ Such a function is described by a sequence of maximal optimal sets that change according to $c$ value.
The corner points depend on the sequence of such solutions.

The algorithm iteratively constructs $PPTP(N^{(k+1)},c)$ by using
$PPTP(N^{(k)},c)$. Starting point is the function $PPTP(N^{(1)},c)$.
It is important to notice that, at each iteration, the description of $PPTP(N^{(k)},c)$ contains at most $k+1$ solutions.
We use  $D^{(k)}$ to indicate the description of function $PPTP(N^{(k)},c)$.
In particular, we indicate as $D^{(k)}.S$  and $D^{(k)}.S_i$
the list of maximal optimal solutions and the
$i$-th solution in such a list, respectively.
Finally,  $D^{(k)}.C_{i}^{\min}$ and $D^{(k)}.C_{i}^{\max}$ are the  minimum and maximum values of  $c$ for which
$D^{(k)}.S_{i}$ is optimal for the problem on $N^{(k)}$.

\bigskip

It follows that:

\begin{itemize}
	\item $D^{(k)}.S_{1}=N^{(k)}$ (see Proposition \ref{PROP_POP_C}
	point \ref{PROP_POP_C_1})
	
	\item $D^{(k)}.S_{|D^{(k)}.S|}=\emptyset$ (see Proposition \ref{PROP_POP_C} point
	\ref{PROP_POP_C_2})
	
	\item $D^{(k)}.C^{\min}_i =
	\left\{
	\begin{array}
	[l]{ll}%
	0 & \text{ for } i=1\\
	c\text{ such that } G(D^{(k)}.S_{i},c)=G(D^{(k)}.S_{i-1},c) & \text{ for } i=2,\ldots,|D^{(k)}.S|\\
	\end{array}
	\right.
	$
	
	\item $D^{(k)}.C^{\max}_i =
	\left\{
	\begin{array}
	[l]{ll}%
	c\text{ such that } G(D^{(k)}.S_{i},c)=G(D^{(k)}.S_{i+1},c) & \text{ for } i=1,\ldots,|D^{(k)}.S|-1\\
	+\infty & \text{ for } i=|D^{(k)}.S|
	\end{array}
	\right.
	$
	
\end{itemize}

\bigskip

{
	\begin{algorithm}[!htbp]
		\caption{\textsc{PPTPLine}($N$)}\label{alg:PopLine}
		\begin{algorithmic}[1]
			\State $\text{Set } h= |N| \text{ and  } D^{(0)}.S=[\emptyset], D^{(0)}.C^{\min}_1=0, D^{(0)}.C^{\max}_1=\infty$
			\For{$k = 0$  \text{\bf  to } $h-1$}
			\State $E^{(k)} \gets \{D^{(k)}.S_{i}\cup\{v_{k+1}\}\ |\ i=1,...,|D^{(k)}.S|\}$
			\State $D^{(k+1)} \gets \textsc{JointSortFilter}(D^{(k)},E^{(k)})$
			\label{PopLine:iterative}
			\EndFor
			\State \Return $(D^{(n)})$ \label{kernelp:returnsol}
		\end{algorithmic}
	\end{algorithm}
}

Function  $\textsc{JointSortFilter}$ executes  the main task of creating the upper envelope
of the solutions in $D^{(k)}$ and $E^{(k)}$.
More precisely, it consider a set  $O(k)$ of solutions: the ones coming from  $D^{(k)}$ are already sorted and define a base of the envelope for $c\in\lbrack0,+\infty)$;\ the ones coming from $E^{(k)}$ are added one at a time to the existing envelope by modifying it accordingly.
All solutions of  $E^{(k)}$ that do not modify the envelope are discarded along with the ones of the existing envelope
that are dominated by the insertion of a new solution.

\begin{theorem}
	Algorithm \textsc{PPTPLine}$(N)$ provides the description of  function $POP(N,c)$ in
	a computational time $O(n^{3}).$
\end{theorem}

\begin{proof}
	Initially, the algorithm provides the description of
	$POP(N^{(1)},c).$ By construction there exists exactly two optimal maximal solutions: the one that includes the unique customer in $N^{(1)}$ and the void solution.
	The moving from the  (minimal) description $D^{(k)}$ of $PPTP(N^{(k)},c)$
	at the beginning of each iteration, to the  (minimal) description $D^{(k+1)}$ of
	$PPTP(N^{(k+1)},c)$ at the end of the iteration is guaranteed by  Proposition
	\ref{PROP_ESTENSIONE_K_K1}.

	As far as computational complexity is concerned,
	we observe that, since
	$|D^{(k)}.S|$ is $O(k),$ each one of the  $O(k)$ solutions  in $E^{(k)}$ have to be compared with the  $O(k)$ solutions  of the current envelope, which means a computational complexity of
	$O(k^{2})$. Since the operation has to be repeated for
	$k=1,...,n-1$ times, the complexity
	$O(n^{3})$ immediately follows. \end{proof}

\section{Conclusions}\label{sec:conclusion}
In this paper, we analyze the Probabilistic Orienteering Problem for the special case where customers are located on a line.
A straightforward algorithm is devised that  allows to determine the upper envelope of the function describing the problem.
The algorithm takes a polynomial time to find the optimal solution.
As future work, the extension of the problem to more general cases
will be taken into account.

\section*{Appendix}

\begin{proof}[Proof of Proposition \ref{PROP:RICORSIONE_1}]
	
	\begin{itemize}
		\item Proof of formula (\ref{RICORSIONE_R})
		\begin{align*}
		R(S)=\sum_{v_{i}\in S}\pi_{i}p_{i} =\left(\sum_{v_{i}\in S^{\prime}}\pi_{i}p_{i}\right) + \pi^{S}p^{S} = R(S^{\prime})+\pi^{S}p^{S};
		\end{align*}
		
		\item Proof of formula (\ref{RICORSIONE_L})
		\begin{align*}
		L(S)= & \sum_{v_{i}\in S}\left[  l_{i}\pi_{i}\cdot{\displaystyle\prod_{{v_j}\in S,\ j>i}}(1-\pi_{j})\right]\\%
		= & \sum_{v_{i}\in S^{\prime}}\left[  l_{i}\pi_{i}\cdot{\displaystyle\prod_{{v_j}\in S,\ j>i}}(1-\pi_{j})\right]+ l^{S}\pi^{S}\cdot 1\\%
		= & (1-\pi^{S})\cdot\sum_{v_{i}\in S^{\prime}}\left[  l_{i}\pi_{i}\cdot{\displaystyle\prod_{{v_j}\in S^{\prime},\ j>i}}(1-\pi_{j})\right]+ l^{S}\pi^{S}\cdot 1\\%
		= & (1-\pi^{S})\cdot L(S^{\prime})+l^{S}\pi^{S};
		\end{align*}
		
		\item Proof of formula (\ref{RICORSIONE_C})
		\begin{align*}
		C(S,c)  &  =cL(S)\\
		&  =c(1-\pi^{S})L(S^{\prime})+\pi^{S}l^{S}c\\
		&  =C(S^{\prime},c(1-\pi^{S}))+\pi^{S}l^{S}c;
		\end{align*}
		
		\item Proof of formula (\ref{RICORSIONE_G})
		\begin{align*}
		G(S,c)  &  =R(S)-C(S,c)\\
		&  =R(S^{\prime})+\pi^{S}p^{S}-(C(S^{\prime},c(1-\pi^{S}))+\pi^{S}l^{S}c)\\
		&  =\left[  R(S^{\prime})-(1-\pi^{S})cL(S^{\prime})\right]  +\pi^{S}%
		(p^{S}-l^{S}c)\\
		&  =G(S^{\prime},(1-\pi^{S})c)+\pi^{S}(p^{S}-l^{S}c).
		\end{align*}
	\end{itemize}
\end{proof}

\bigskip
\begin{proof}[Proof of Proposition \ref{PROP_SUBSET_LENGTH_INEQUALITY}]
	Let $S_1$ and $S_2$ be two set of customers such that $S_1 \subseteq S_2$. we show that $L(S_1)\le L(S_2)$ by induction on cardinality of $S_2$ using equality (\ref{RICORSIONE_L}).
	
	\textbf{\textit{Base cases}}. If $|S_2|=0$, both $S_1$ and $S_2$ are empty and we get $L(S_1)=0\le L(S_2)=0$. If $|S_2|=1$, then either $S_1=\emptyset$ or $S_1=S_2$; and we get $L(S_1)=0\le L(S_2);$ and $L(S_1)=L(S_2)$, respectively.
	
	\textbf{\textit{Induction hypothesis}}. Let us assume that the property holds for all $|S_2|\leq n$ and show that it hold also for $|S_2|=n+1$.
	
	\textbf{\textit{Induction step}}. Let $n>1$ be the cardinality of $S_2$. Then either $v^{S_2}\in S_1$ or $v^{S_2} \notin S_1$.
	
	If $v^{S_2}\in S_1$, we have $v^{S_2}=v^{S_1}=v^{S}$, and $L(S_1)=(1-\pi^{S}\cdot L(S_1^{\prime})+ l^{S}\pi^{S}\le L(S_2^{\prime})+ l^{S}\pi^{S}=L(S_2)$; the last inequality comes from induction hypothesis because $|S_2^{\prime}|=|S_2|-1=n$.
	
	If $v^{S_2} \notin S_1$, we have $S_1\subseteq S_2$ and $L(S_1)\leq L(S_2^{\prime})\leq L(S_2^{\prime})+l^{S_2}\pi^{S_2}=L(S_2)$; the first inequality comes from induction hypothesis because $|S_2^{\prime}|=|S_2|-1=n$.
\end{proof}

\bigskip
\begin{proof}[Proof of Proposition \ref{PROP_INSIEMI}]
	\begin{itemize}
		\item Proof of formula \eqref{PROP_INSIEMI_R}
		Equality follows straightforwardly from the following decompositions:
		\begin{align*}
		R(S_1)  &  =\sum_{v_{i}\in S_1\backslash S_2}\pi_{i}p_{i}+\sum_{v_{i}\in S_1\cap S_2}%
		\pi_{i}p_{i},\\
		R(S_2)  &  =\sum_{v_{i}\in S_2\backslash S_1}\pi_{i}p_{i}+\sum_{v_{i}\in S_1\cap S_2}%
		\pi_{i}p_{i},\\
		R(S_1\cup S_2)  &  =\sum_{v_{i}\in S_1\backslash S_2}\pi_{i}p_{i}+\sum_{v_{i}\in
			S_2\backslash S_1}\pi_{i}p_{i}+\sum_{v_{i}\in S_1\cap S_2}\pi_{i}p_{i},\\
		R(S_1\cap S_2)  &  =\sum_{v_{i}\in S_1\cap S_2}\pi_{i}p_{i};
		\end{align*}

		\item Proof of formula (\ref{PROP_INSIEMI_L})
		We first discuss two special cases.
		\begin{enumerate}
			\item[a)]
			When $S_2\subseteq S_1$ (or $S_1\subseteq S_2$)
			the property is trivially true since $S_1\cup S_2=S_1$ and $S_1\cap S_2=S_2$ (or $S_1\cup S_2=S_2$ and $S_1\cap S_2=S_1$) and equality boils down to an identity.

			\item[b)]
			When $S_1\cap
			S_2=\emptyset$ the property boils down to the inequality  $L(S_1\cup S_2)\leq
			L(S_1)+L(S_2)$ which comes from the comparison of the following expressions:
			\begin{align*}
			L(S_1)  &  =\sum_{v_{i}\in S_1}\pi_{i}l_{i}%
			{\displaystyle\prod_{{j}\in S_1,\ j>i}}
			(1-\pi_{j}),\\
			L(S_2)  &  =\sum_{v_{i}\in S_2}\pi_{i}l_{i}%
			{\displaystyle\prod_{{j}\in S_2,\ j>i}}
			(1-\pi_{j})\\
			L(S_1\cup S_2)  &  =\sum_{v_{i}\in S_1}\pi_{i}l_{i}%
			{\displaystyle\prod_{{j}\in S_1\cup S_2,\ j>i}}
			(1-\pi_{j})+\sum_{v_{i}\in S_2}\pi_{i}l_{i}%
			{\displaystyle\prod_{{j}\in S_1\cup S_2,\ j>i}}
			(1-\pi_{j}).
			\end{align*}
			where inequalities
			\begin{align*}
			\sum_{v_{i}\in S_1}\pi_{i}l_{i}{\displaystyle\prod_{{j}\in S_1\cup S_2,\ j>i}}(1-\pi_{j})  &  \leq\sum_{v_{i}\in S_1}\pi_{i}l_{i}{\displaystyle\prod_{{j}\in S_1,\ j>i}}(1-\pi_{j}),\\
			\sum_{v_{i}\in S_2}\pi_{i}l_{i}{\displaystyle\prod_{{j}\in S_1\cup S_2,\ j>i}}(1-\pi_{j}) &  \leq\sum_{v_{i}\in S_2}\pi_{i}l_{i}{\displaystyle\prod_{{j}\in S_2,\ j>i}}(1-\pi_{j})
			\end{align*}
			hold because each term in the left summations has more factor less than $1$.
		\end{enumerate}
		
		\medskip
		We conduct the proof by induction on the cardinality of $S=S_1\cup S_2$ and make use of notation of equation \eqref{RICORSIONE_L} which we recall here for reader's convenience
		$$L(S)=(1-\pi^{S})L(S^{\prime})+\pi^{S} l^{S};$$
		where we indicate as  $v^{S}$ the farthest customer of a given set $S$ (with location  $l^{S},$ profit $p^{S}$ and  probability $\pi^{S}$), and define    $S^{\prime}=S\backslash\{v^{S}\}$.
		
		\textbf{\textit{Base cases}}.
		If $|S|\leq2$ we have (w.l.o.g.) $S_2\subseteq S_1$ and/or $S_1\cap S_2=\emptyset$. Thus, we fall in one of cases a,b discussed above.
		
		\textbf{\textit{Induction hypothesis}}.
		We assume that the property holds for $|S|=|S_1\cup S_2|\leq n$ and show that it holds also for $|S|=n+1$.
		
		\textbf{\textit{Induction step}}.
		Given $|S_1\cup S_2|=n+1$, we assume w.l.o.g. that $v^{S}\in S_1$. The following two cases may occur:
		
		\begin{itemize}
			\item[a)] $v^{S}\in S_2$ (thus $v^{S}\in S_1\cap S_2$).
			From \eqref{RICORSIONE_L}, we get:
			\begin{align*}
			L(S_1\cup S_2)  &  =\pi^{S} l^{S}+(1-\pi^{S})L(S_1^{\prime}\cup S_2^{\prime})\\
			L(S_1\cap S_2)  &  =\pi^{S} l^{S}+(1-\pi^{S})L(S_1^{\prime}\cap S_2^{\prime})\\
			L(S_1)  &  =\pi^{S} l^{S}+(1-\pi^{S})L(S_1^{\prime})\\
			L(S_2)  &  =\pi^{S} l^{S}+(1-\pi^{S})L(S_2^{\prime})
			\end{align*}
			and thus
			$L(S_1\cup S_2)+L(S_1\cap S_2)=2\pi^{S} l^{S}+(1-\pi^{S})\left[  L(S_1^{\prime}\cup S_2^{\prime
			})+L(S_1^{\prime}\cap S_2^{\prime})\right]
			$,
			where $|S_1^{\prime}\cup S_2^{\prime}|=n.$
			Then, by induction hypothesis: 
			\begin{align*}
			L(S_1\cup S_2)+L(S_1\cap S_2)  &  \leq2\pi^{S} l^{S}+(1-\pi^{S})\left[  L(S_1^{\prime})+L(S_2^{\prime
			})\right] \\
			&  \leq\left(  \pi^{S} l^{S}+(1-\pi^{S})L(S_1^{\prime})\right)  +\left(  \pi^{S} l^{S}+(1-\pi^{S}
			)L(S_2^{\prime})\right) = L(S_1)+L(S_2).
			\end{align*}

			\item[b)] $v^{S}\notin S_2$, thus
			\[
			v^{S}\notin S_1 \cap S_2 \ \text{with }
			S^{\prime}=S_1^{\prime}\cup S_2 \ \text{and }
			|S^{\prime}| = n.
			\]
			From \eqref{RICORSIONE_L}, we get $L(S_1) =\pi^{S} l^{S}+(1-\pi^{S})L(S_1^{\prime})$,
			and then
						\[
			L(S_1^{\prime})=L(S_1)+\pi^{S}\left(  L(S_1^{\prime})-l^{S}\right).
			\]
			By using the induction hypothesis  ($|S_1^{\prime}\cup S_2|=n$) we get
			\begin{align*}
			L(S_1^{\prime}\cup S_2)+L(S_1^{\prime}\cap S_2)  &  \leq L(S_1^{\prime})+L(S_2)%
			= L(S_1)+\pi^{S}\left(  L(S_1^{\prime})-l^{S}\right)  +L(S_2)
			\end{align*}
			Now, the final result is obtained by substituting back and recalling that $L(S_1\cap
			S_2)=L(S_1^{\prime}\cap S_2)$ (from $v^{S}\notin S_2$) and  $L(S_1^{\prime})\leq L(S_1^{\prime}\cup S_2)$ (from Proposition \ref{PROP_SUBSET_LENGTH_INEQUALITY}):%
			\begin{align*}
			L(S_1\cup S_2)+L(S_1\cap S_2)  &  =\pi^{S} l^{S}+(1-\pi^{S})L(S_1^{\prime}\cup S_2)+L(S_1^{\prime}\cap
			S_2)\\
			&  =\pi^{S} l^{S}-\pi^{S} L(S_1^{\prime}\cup S_2)+L(S_1^{\prime}\cup S_2)+L(S_1^{\prime}\cap S_2)\\
			&  \leq\pi^{S} l^{S}-\pi^{S} L(S_1^{\prime}\cup S_2)+L(S_1^{\prime})+L(S_2)\\
			&  =\pi^{S} l^{S}-\pi^{S} L(S_1^{\prime}\cup S_2)+L(S_1)+\pi^{S}\left(  L(S_1^{\prime})-l^{S}\right)
			+L(S_2)\\
			&  =\pi^{S}\left(  l^{S}-L(S_1^{\prime}\cup S_2)+L(S_1^{\prime})-l^{S}\right)  +L(S_1)+L(S_2)\\
			&  =\pi^{S}\left(  L(S_1^{\prime})-L(S_1^{\prime}\cup S_2)\right)  +L(S_1)+L(S_2)\\
			&  \leq L(S_1)+L(S_2).
			\end{align*}
		\end{itemize}
		The same considerations hold for the case  $v\in S_2\backslash S_1$.
	\end{itemize}
\end{proof}

\end{document}